\documentclass{article}
\usepackage{amsmath,amsfonts,amssymb}

\begin{document}
\title{Limit Representations of Riemann's Zeta Function}

\author{Djurdje Cvijovi\'{c} and Hari  M.  Srivastava}
\date{}
\maketitle

\begin{abstract}
In this paper   it is shown that Riemann's zeta function $\zeta(s)$ admits two limit representations when $\Re{(s)}>1.$  Each of these limit representations is deduced by using simple arguments  based upon the classical Tannery's (limiting) theorem for series.
\end{abstract}

\vskip5mm
\section{Introduction.}
Riemann's zeta function $\zeta(s)$ is  a complex-valued function of a complex variable $s$ and is holomorphic everywhere in the complex $s$-plane except at the point $s = 1$ where a first-order pole exists with residue equal to 1. It is, as usual, defined as the analytic continuation of the function given by the sum of the following  series:

\begin{equation}
\zeta(s):=\left\{ \aligned  &\sum_{n=1}^\infty
\frac{1}{n^s}=\frac{1}{1-2^{-s}} \sum_{n=1}^\infty
   \frac{1}{(2n-1)^s} \qquad \hskip 2 mm  \big(\Re(s)>1\big) \\
   & \frac{1}{1-2^{1-s}} \, \sum_{n=1}^\infty \frac{(-1)^{n-1}}{n^s}
      \qquad   \hskip 18 mm \big(\Re(s)>0;\; s \not= 1\big).
\endaligned \right.
\end{equation}

\noindent Moreover, a number of other infinite series, infinite products, improper integrals, complex contour integrals and closed-form expressions based upon the Euler-Maclaurin summation formula may be used to represent $\zeta(s)$ in certain regions of  the complex $s$-plane (see, for details, \cite{Apostol0} and \cite{SC}). Some illustrative examples are given below.

Euler's product formula for the zeta function:

\begin{equation}
\zeta(s)  = \prod_{m\,=1}^{\infty}\frac{1}{1-p_m^{-s}} \qquad \big(\Re(s) > 1\big),
\end{equation}

\noindent where $p_m$ is $m$th prime number;
\begin{equation}
\zeta(s) = \frac{1}{\Gamma(s)}\int_{0}^{\infty}\frac{x^{s-1}} {e^x-1}  \,dx\qquad \big(\Re(s) > 1\big);
\end{equation}

\begin{equation}
\zeta(s) = \frac{1}{s-1} + \sum_{n\,=0}^{\infty}(-1)^n \,\frac{\gamma_n}{n!} \,(s-1)^n \qquad(s\neq 1),
\end{equation}

\noindent where
\begin{equation}\gamma_n = \lim _{m \rightarrow \infty }\left( \sum_{k\,=1}^{m} \frac{(\log k)^n}{k} - \frac{(\log m)^{n+1}}{n+1}\right)\quad\big(n \in \mathbb{N}:=\{1,2,3,\cdots\} \big);
\end{equation}

\begin{equation}
\zeta(s)= \sum_{m\,=1}^n \frac{1}{m^s} + \frac{n^{1-s}}{s-1} -s \int_{n}^{\infty} \frac{x- \lfloor x\rfloor}{x^{s+1}}\,dx\qquad\big(\Re(s) > 0; n \in \mathbb{N}\big),
\end{equation}

\noindent where $\lfloor x\rfloor$ stands for the floor function which gives the largest integer less than or equal to $x\in \mathbb{R}$;

\begin{equation}
\zeta(s) = \frac{\Gamma(1-s)}{2 \pi \imath}\int_{-\infty}^{(0+)}\frac{z^{s-1}} {e^{-z}-1}  \,dz
\qquad (s\in \mathbb{C}\setminus \mathbb{N}),
\end{equation}

\noindent where the contour of integration is a loop around the negative real axis; it starts at $-\infty$, encircles the origin once in the positive (counter-clockwise) direction without enclosing any of the points
$$z = \pm \, 2 \pi \imath, \pm \,4 \pi \imath,\cdots,$$
and returns to $ -\infty$.

All of the above and many other representations of $\zeta(s)$ are known for a considerable time. For an exhaustive list of such and other representations of $\zeta(s)$, the interested reader is referred (for instance) to \cite{Apostol0} and \cite{SC}. This is not surprising since there is a long and rich history of research on Riemann's zeta function $\zeta(s)$ that goes back to Euler in 1735 (see, for details, \cite{Ayoub} and \cite{Sri2007-JKMS}). What is surprising, however, is that it has not been noticed hitherto that $\zeta(s)$ admits two limit representations which are asserted here by the following theorem.

\vskip 2mm\noindent
{\bf Theorem.} {\it Suppose that $s$ is a complex number and let $m,$ $n,$ $p$ and $q$ be  nonnegative integers. Then$,$ for $ \Re(s) > 1,$  the values of Riemann's zeta function $\zeta(s)$ are given by}
\begin{align}
&\textup{(a)} \quad \zeta(s) = \lim _{q \rightarrow \infty } \left(\frac{\pi}{2\, q + m}\right)^{s} \sum_{p\,= 1}^{\lfloor (2\, q  +  n  -  1)/2\rfloor} \cot^{s} \left(\frac{p\, \pi}{2\, q + n}\right)\hskip 22mm
\\
& \hskip -3mm {\it and}\nonumber
\\
&\textup{(b)} \quad \zeta(s) =  \lim _{q \rightarrow \infty } \left(\frac{\pi}{2\, q + m}\right)^{s} \sum_{p\,= 1}^{\lfloor (2\, q  +  n  -  1)/2\rfloor} \csc^{s} \left(\frac{p\, \pi}{2\, q + n}\right)
\end{align}
$$(p\in\mathbb{N},  m\in \mathbb{N}_0:=\mathbb{N}\cup \{0\};\; n = 0 \quad \textit{and}\quad q \in \mathbb{N}\setminus \{1\};\; n \in \mathbb{N} \quad \textit{and}\quad q \in \mathbb{N}).$$

\noindent{\bf Remark 1.}
Several special cases of the limit relationships (8) and (9) involving $\zeta(2 n)$ and $\zeta(2 n +1)$ when $n\in \mathbb{N}$ can be found in the literature. For example, by making use of elementary arguments and complex function theory, respectively, the following two limit formulas were established by Williams \cite[p. 273, Lemma 1; p. 275, Lemma 2]{Williams}:

\begin{equation}\zeta(2n)  = \lim _{q \rightarrow \infty } \left(\frac{\pi}{2 q}\right)^{2 n}
\sum_{p\,= 1}^q \cot^{2 n} \left(\frac{p\, \pi}{2 q+1}\right)
\end{equation}

\noindent and

\begin{equation}
\zeta(2n)  = \lim _{q \rightarrow \infty } \left(\frac{\pi}{2 q + 1}\right)^{2 n}
\sum_{p\,= 1}^q \cot^{2 n} \left(\frac{p\, \pi}{2 q+1}\right).
\end{equation}

\noindent On the other hand, an elementary proof for the following limit formula was given by Apostol \cite[p. 430, Eq. (16)]{Apostol}:

\begin{align}\zeta(2n +1) = \lim _{q \rightarrow \infty } \left(\frac{\pi}{2 q}\right)^{2 n +1}
\sum_{p\,= 1}^q \cot^{2 n+1} \left(\frac{p\, \pi}{2 q+1}\right),
\end{align}

\noindent who also found an asymptotic expansion of the finite sum in (10) \cite[p. 428, Eq. (7)]{Apostol}, which readily leads to Euler's celebrated relation:

\begin{equation}
\zeta(2 n) = (-1)^{n+1}\; \frac{(2 \pi)^{2 n}}{2 \cdot (2 n)!}\, B_{2 n} \qquad (n \in \mathbb{N}_0)
\end{equation}

\noindent between the even-indexed Bernoulli numbers $B_{2 n}$ and the values of the $\zeta(2n).$  It should be noted that (12) can be proven by following Williams' arguments used in the case of (10) with necessary changes.
Cvijovi\'{c} {\it et al.} \cite{Cvijovic} resorted to the calculus of residues in order to derive the cotangent finite sum in (10) and some other related sums in a closed form. Furthermore, as immediate consequences of their results, Williams' limit formulas (10) and (11) as well as the following three related limit formulas for $\zeta(2n)$ were obtained by Cvijovi\'{c} {\it et al.} \cite[p. 206, Theorem 2]{Cvijovic}:

\vskip2mm
\begin{equation}
\zeta(2n) = \lim _{q \rightarrow \infty } \left(\frac{\pi}{2 q}\right)^{2n}
\sum_{p\,= 1}^{q-1} \cot^{2n} \left(\frac{p\, \pi}{2q}\right),
\end{equation}
\begin{equation}
\zeta(2n)= \lim _{q \rightarrow \infty } \left(\frac{\pi}{2 q}\right)^{2n}
\sum_{p\,= 1}^{q-1} \csc^{2n} \left(\frac{p\, \pi}{2q}\right)
\end{equation}
and
\begin{equation}
\zeta(2n)=\lim _{q \rightarrow \infty } \left(\frac{\pi}{2 q}\right)^{2 n}
\sum_{p\,= 1}^q \csc^{2 n} \left(\frac{p\, \pi}{2 q+1}\right).
\end{equation}

\section{Demonstration of the Theorem.}

After numerous unsuccessful attempts to generalize the limit formulas (10) to (12) and (14) to (16), we have encountered an old and almost elementary result of the classical analysis, known as Tannery's (limiting) theorem for series (see \cite[p. 292]{Tannery}, \cite[pp, 123 and 124]{Bromwich}, \cite[pp. 371 and 372] {MacRobert} and \cite[pp. 199 and 200]{Hofbauer}), which indeed provides a simple and direct proof of our Theorem. We first state Tannery's theorem here {\it without} proof, noting that its standard application is to show that the following two usual definitions of $e^x$ are the same:

\begin{equation}
\lim _{n\rightarrow \infty } \left(1 + \frac{x}{n}\right)^n = \lim _{n\rightarrow \infty } \sum_{k\,=0}^n \binom{n}{k} \left(\frac{x}{n}\right)^{k} = \sum_{k\,= 0}^{\infty} \frac{x^k}{k!}= :e^x.
\end{equation}

\vskip2mm\noindent
{\bf Tannery's Theorem} ({\it cf.} \cite[pp. 371 and 372]{MacRobert}). {\it For a given double sequence} $\{f_m(n)\}_{m,n\in \mathbb{N}_0},$ {\it suppose that each of the following two conditions is satisfied for any fixed} $m \in \mathbb{N}_0:$

\vskip 3mm

\textup(i)   $\lim\limits _{n \rightarrow \infty } f_m(n) = f_m ;$

\vskip 2mm
\textup(ii) $\left|f_m(n)\right|\leqq M_m,$ {\it where $M_m > 0$ is independent of}  $\,n\,$ {\it and the infinite series$:$} $$\sum_{m\,=0}^{\infty} M_m$$ {\it is convergent.

Then the following limit relationship holds true$:$}

\begin{equation}
\lim _{n \rightarrow \infty}\; \sum_{m\,= 0}^{\alpha(n)} f_m (n)= \sum_{m\,= 0}^{\infty} f_m,
\end{equation}
{\it where} $\{\alpha(n)\}_{n\in \mathbb{N}_0}$ {\it is a monotonically increasing integer-valued sequence which tends to infinity as} $n \rightarrow \infty$.

\vskip 2mm In what follows, it is assumed that $m$, $n$, $p$ and $q$ are nonnegative integers and we also {\em suppose, for a moment,  that $s$ is a real number}.

To prove Part (a) when $s > 1$, we consider the following double sequence which appears in (8):
\begin{equation}
\Phi_{p} (q \,|\, m, n, s) :=  \left(\frac{\pi}{2 \,q + m}\right)^s \cot^s \left( \frac{p\, \pi}{2\, q + n}\right),
\end{equation}
where $p$ and $q$ are indices and  $m$, $n$ and $s$ are parameters with
\begin{equation}
p = 1, \cdots, \Big\lfloor\frac{2\, q + n - 1}{2}\Big\rfloor
\end{equation}
$$(p\in\mathbb{N}, m\in \mathbb{N}_0;\; n = 0 \quad \textit{and}\quad q \in \mathbb{N}\setminus \{1\};\; n \in \mathbb{N} \quad \textit{and}\quad q \in \mathbb{N}).$$
We now show that the double sequence defined by (19) satisfies  the  conditions (i) and (ii) of Tannery's theorem.

\vskip 2mm\noindent
{\bf Condition (i) of Tannery's Theorem.} For a fixed $p$ and fixed $m,$  $n$ and $s$, in view of the well-known facts that
$$\tan x \sim x \quad (x \rightarrow 0),\;\; \frac{p\, \pi}{2\, q + n} \rightarrow 0 \quad \text{and} \quad
\frac{2 \,q + n}{2\, q + m}  \rightarrow 1 \quad (q \rightarrow \infty),$$
it follows without difficulty that
\begin{align}
\lim_{q \rightarrow \infty } \left[\frac{\pi}{2\, q + m}\, \cot \left( \frac{p\, \pi}{2\, q + n}\right)\right]^s  = \left[ \lim _{q \rightarrow \infty } \left(\frac{\pi}{2 \,q + m}\right)\left( \frac{2 \,q + n}{p \,\pi}\right)\right]^s = \frac{1}{p^s}.
\end{align}

\vskip 2mm\noindent
{\bf Condition (ii) of Tamnnery's Theorem.} In this case, we recall the following elementary inequality \cite[p. 75, Entry 4.3.80]{Abramowitz}: $$\sin x< x<\tan x \qquad \left(0 <x < \frac{\pi}{2}\right),$$
which yields $$0 < \cot x < \frac{1}{x} \qquad \left(0 <x < \frac{\pi}{2}\right),$$
so that, since [{\it cf.} Equation (20)]
$$0 < \frac{p \pi}{2 \,q + n} < \frac{\pi}{2},$$
we have
\[0 < \cot \left(\frac{p \,\pi}{2 \,q + n}\right)< \frac{2 \,q + n}{p \,\pi},
\]
which, upon noticing that $\pi/(2 \,q + m) > 0$, becomes
\begin{equation}
0 < \frac{\pi}{2 \,q + m} \,\cot \left(\frac{p\, \pi}{2\, q + n}\right)< \frac{2\, q + n}{2\, q + m}\cdot \frac{1}{p}.
\end{equation}
Moreover, from the graph of the function defined by $$f(x):=\frac{2 x +n}{2 x+m},$$
it is easily concluded that
\begin{equation}
\frac{2 q +n}{2 q+m}\leqq C_{m,n}:=
\begin{cases}
1  & \quad (n\leqq m) \\
\\
\frac{1+n}{1+m} & \quad (n > m)
\end{cases}
\qquad(m,n\in \mathbb{N}_0;\; q\in \mathbb{N}).
\end{equation}
Thus, by making use of these last two equations (22) and (23), we find for $s > 0$ that
\begin{equation}\left| \left(\frac{\pi}{2 \,q + m}\right)^s \cot^s \left( \frac{p\, \pi}{2 \,q + n}\right)\right|  = \left[\frac{\pi}{2\, q + m}\, \cot \left(\frac{p \pi}{2\, q + n}\right)\right]^s \leqq C_{m,n}^s   \frac{1}{p^s}.
\end{equation}
Clearly, therefore, the condition (ii) of Tannery's theorem is fulfilled only when $ s > 1$, because the infinite series
$$\sum_{p\,= 1}^{\infty}\frac{1}{p^{s}}$$
is then convergent.

In conclusion, we may apply Tannery's theorem to the double sequence
$$\Phi_{p} (q \,|\, m, n, s) \qquad (s>1),$$
given by the equations (19) and (20), since the needed conditions are satisfied and  $\lfloor(2\, q + n -1)/2\rfloor$ is evidently an increasing integer-valued function which tends to infinity as $q\rightarrow \infty$.  Thus, for $ s > 1$, the desired limit formula in (8) follows in view of the limit relationship:

\begin{equation}
\lim _{q \rightarrow \infty } \sum_{p\,= 1}^{\lfloor(2\, q + n - 1)/2\rfloor} \left(\frac{\pi}{2 \,q + m}\right)^s \cot^s \left( \frac{p\, \pi}{2 \,q + n}\right)= \sum_{p\,= 1}^{\infty} \frac{1}{p^s} = \zeta(s)\quad(s > 1).
\end{equation}

To prove Part (b) when $s >1$, we consider the double sequence defined by

\begin{equation}
\Psi_{p}(q \,|\, m, n, s) = \left(\frac{\pi}{2 \,q + m }\right)^s \csc^s \left( \frac{p\, \pi}{2 \,q + n}\right),
\end{equation}

\noindent together with the restrictions on the integers $p,$ $q,$ $m$ and $n$ given in (20) and proceed along the same lines as in the proof of Theorem (a). In the process, it is necessary to employ the following well-known asymptotic relation:  $$\sin x \sim x \qquad (x \rightarrow 0)$$ as well as the inequality \cite[p. 75, Entry 4.3.79]{Abramowitz}:
\[0 < \csc x < \frac{\pi}{2\,x} \qquad \left(0 < x < \frac{\pi}{2}\right).
\]

\noindent In this way, by Tannery's theorem for series, we obtain the limit relationship:

\begin{equation}\lim _{q \rightarrow \infty } \sum_{p \,= 1}^{\lfloor(2\, q + n - 1)/2\rfloor} \left(\frac{\pi}{2 \,q + m}\right)^s \csc^s \left( \frac{p\, \pi}{2 \,q + n}\right)= \sum_{p\,= 1}^{\infty} \frac{1}{p^s} = \zeta(s)\quad(s > 1),
\end{equation}

\noindent which obviously implies  the limit formula (9) which holds true for $s > 1$.

\vskip 2mm\noindent
{\bf Remark 2.} Observe that the limit formulas (25) and (27) are deduced here on the supposition that $s$ is a real number. However, these limit formulas are valid in the entire half-plane $\Re{(s)}>1$, since they may be extended by applying the principle of analytic continuation on $s$ as far as possible.

\section{Concluding Remarks and Observations.}

By suitably applying the Theorem, the above-presented limit formulas in (10) to (12) and (14) to (16) could be generalized so as to be valid for $\Re(s)>1.$ We thus have

\begin{equation}
\zeta(s)  = \lim _{q \rightarrow \infty } \left(\frac{\pi}{2 q}\right)^{s} \sum_{p\,= 1}^q \cot^{s} \left(\frac{p\, \pi}{2 q+1}\right) \qquad (\Re(s)>1),
\end{equation}
\begin{equation}
\zeta(s) = \lim _{q \rightarrow \infty } \left(\frac{\pi}{2 q + 1}\right)^{s}
\sum_{p\,= 1}^q \cot^{s} \left(\frac{p\, \pi}{2 q+1}\right)\qquad (\Re(s)>1),
\end{equation}
\begin{equation}
\zeta(s) =\lim _{q \rightarrow \infty } \left(\frac{\pi}{2 q}\right)^{s}
\sum_{p\,= 1}^{q-1} \cot^{s} \left(\frac{p\, \pi}{2 q}\right)\qquad (\Re(s)>1),
\end{equation}
\begin{equation}
\zeta(s) = \lim _{q \rightarrow \infty } \left(\frac{\pi}{2 q}\right)^{s}
\sum_{p\,= 1}^q \csc^{s} \left(\frac{p\, \pi}{2 q + 1}\right)\qquad (\Re(s)>1)
\end{equation}
and
\begin{equation}
\zeta(s) =\lim _{q \rightarrow \infty } \left(\frac{\pi}{2 q}\right)^{s}
\sum_{p\,= 1}^{q-1} \csc^{s} \left(\frac{p\, \pi}{2 q}\right)\qquad (\Re(s)>1).
\end{equation}

\vskip2mm\noindent
{\bf Remark 3.} We remark that many elementary and special functions possess limit representations and the rather well-known ones are those of the exponential function in (17) and Euler's limit formula for the gamma function \cite[p. 255, Entry 6.1.2]{Abramowitz}:

\begin{equation}
\Gamma(z)= \lim _{n \rightarrow \infty }\frac{n! \,n^z}{z(z+1) (z+2)\cdots (z+ n)}
\end{equation}
$$\big(z\in \mathbb{C} \setminus \mathbb{Z}^{-}_0 \quad (\mathbb{Z}^{-}_0:=\{0, -1, -2, -3, \cdots\})\big).$$

\vskip4mm
\noindent\textbf{Djurdje U.  Cvijovi\'{c}} graduated  from  the  University of Belgrade,  Serbia,  and  received his Ph.D. from  the University of Cambridge, England, in 1994. His main research interests include global optimization, theory of special functions and elementary number theory.

\noindent\textit{\it Atomic Physics Laboratory$,$ Vin\v{c}a Institute of
Nuclear Sciences, P. O. Box $522$, Belgrade, Republic of Serbia}

\noindent\textit{djurdje@vinca.rs}

\end{document}